\documentclass[12pt]{article}%
\usepackage{amsmath}
\usepackage{amsfonts}
\usepackage{amssymb}
\usepackage{graphicx}%
\begin{document}

\title{Some remarks on the Wu-Sprung potential. Preliminary report}
\author{Diego Dominici \thanks{e-mail: dominicd@newpaltz.edu}\\Department of Mathematics\\State University of New York at New Paltz\\75 S. Manheim Blvd. Suite 9\\New Paltz, NY 12561-2443\\USA\\Phone: (845) 257-2607\\
Fax: (845) 257-3571
}

\date{}

\maketitle

\section{Introduction}

Let $V(x)$ be the solution of the equation%
\begin{equation}
\pi x=\sqrt{V-V_{0}}\ln\left(  \frac{V_{0}}{2\pi e^{2}}\right)  +\sqrt{V}%
\ln\left[  \frac{\sqrt{V}+\sqrt{V-V_{0}}}{\sqrt{V}-\sqrt{V-V_{0}}}\right]  ,
\label{EqV}%
\end{equation}
where $V_{0}>0.$ Rationalizing the denominator, we can write (\ref{EqV}) in the
more convenient form%
\begin{equation}
\pi x=\sqrt{V-V_{0}}\ln\left(  \frac{V_{0}}{2\pi e^{2}}\right)  +\sqrt{V}%
\ln\left[  \frac{\left(  \sqrt{V}+\sqrt{V-V_{0}}\right)  ^{2}}{V_{0}}\right]
. \label{EqV1}%
\end{equation}

Setting
\begin{equation}
V=V_{0}\sec^{2}(\theta),\quad-\frac{\pi}{2}<\theta<\frac{\pi}{2}
\label{EqVtheta}%
\end{equation}
in (\ref{EqV1}), we obtain%
\begin{equation}
\frac{\pi x}{\sqrt{V_{0}}}=\tan\left(  \theta\right)  \ln\left(  \frac{V_{0}%
}{2\pi e^{2}}\right)  +2\sec\left(  \theta\right)  \ln\left[  \sec\left(
\theta\right)  +\tan\left(  \theta\right)  \right]  . \label{EqXtheta}%
\end{equation}
It follows from (\ref{EqXtheta}) that $x=0$ if $\theta=0$ or $\theta
=\theta_{0}$, with%
\begin{equation}
V_{0}=2\pi e^{2}\left[  \sec\left(  \theta_{0}\right)  +\tan\left(  \theta
_{0}\right)  \right]  ^{-2\sec\left(  \theta_{0}\right)  }. \label{Eqtheta0}%
\end{equation}
In the first case $x=0$ and $V=V_{0},$ and in the second $x=0$ and
$V=V_{0}\sec^{2}(\theta_{0})>V_{0}.$

\begin{figure}[ptb]
\begin{center}
\rotatebox{270} {\resizebox{12cm}{!}{\includegraphics{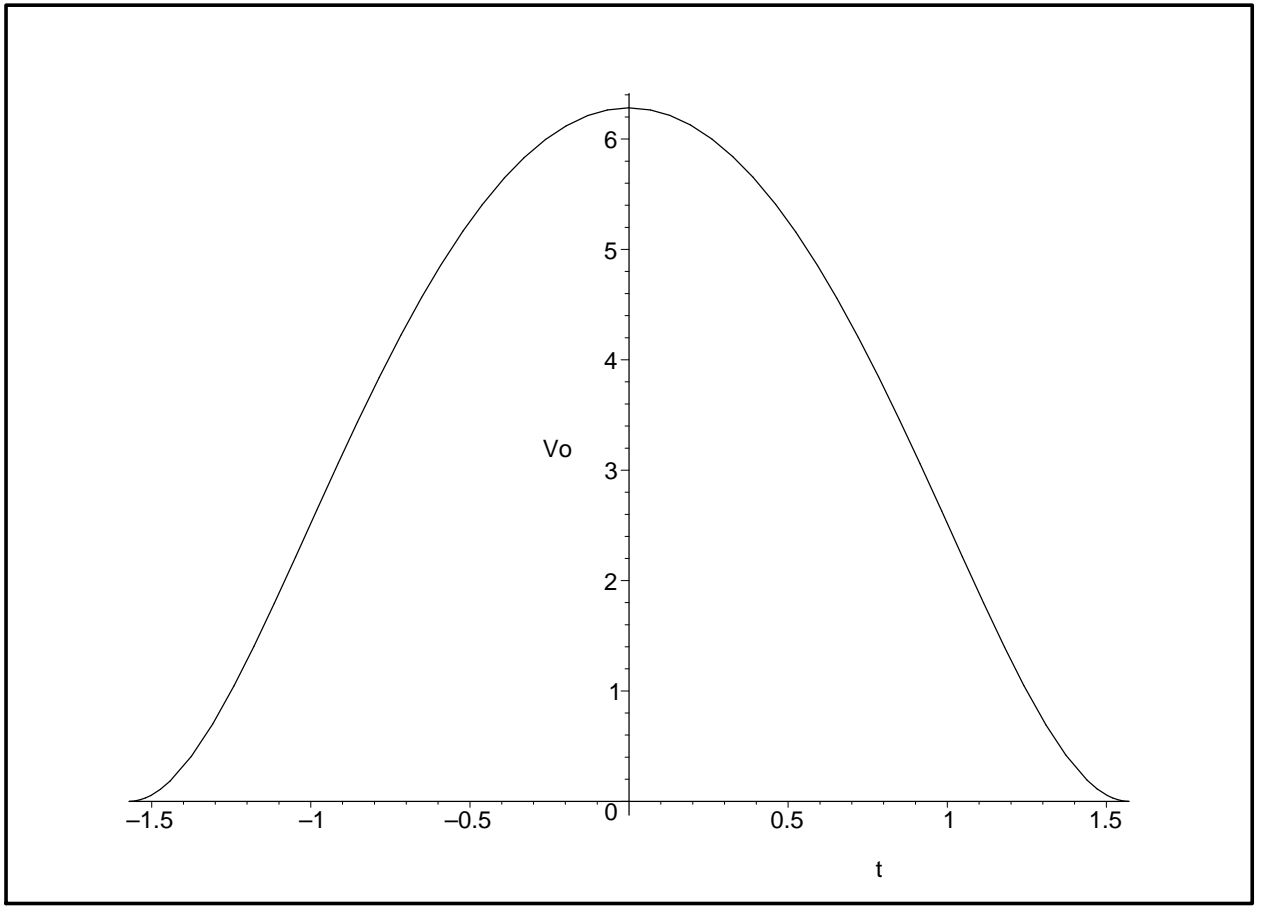}}}
\end{center}
\caption{A sketch of $V_{0}$ as a function of $\theta_{0}$.}%
\label{Vo}%
\end{figure}

\begin{figure}[ptb]
\begin{center}
\rotatebox{270} {\resizebox{12cm}{!}{\includegraphics{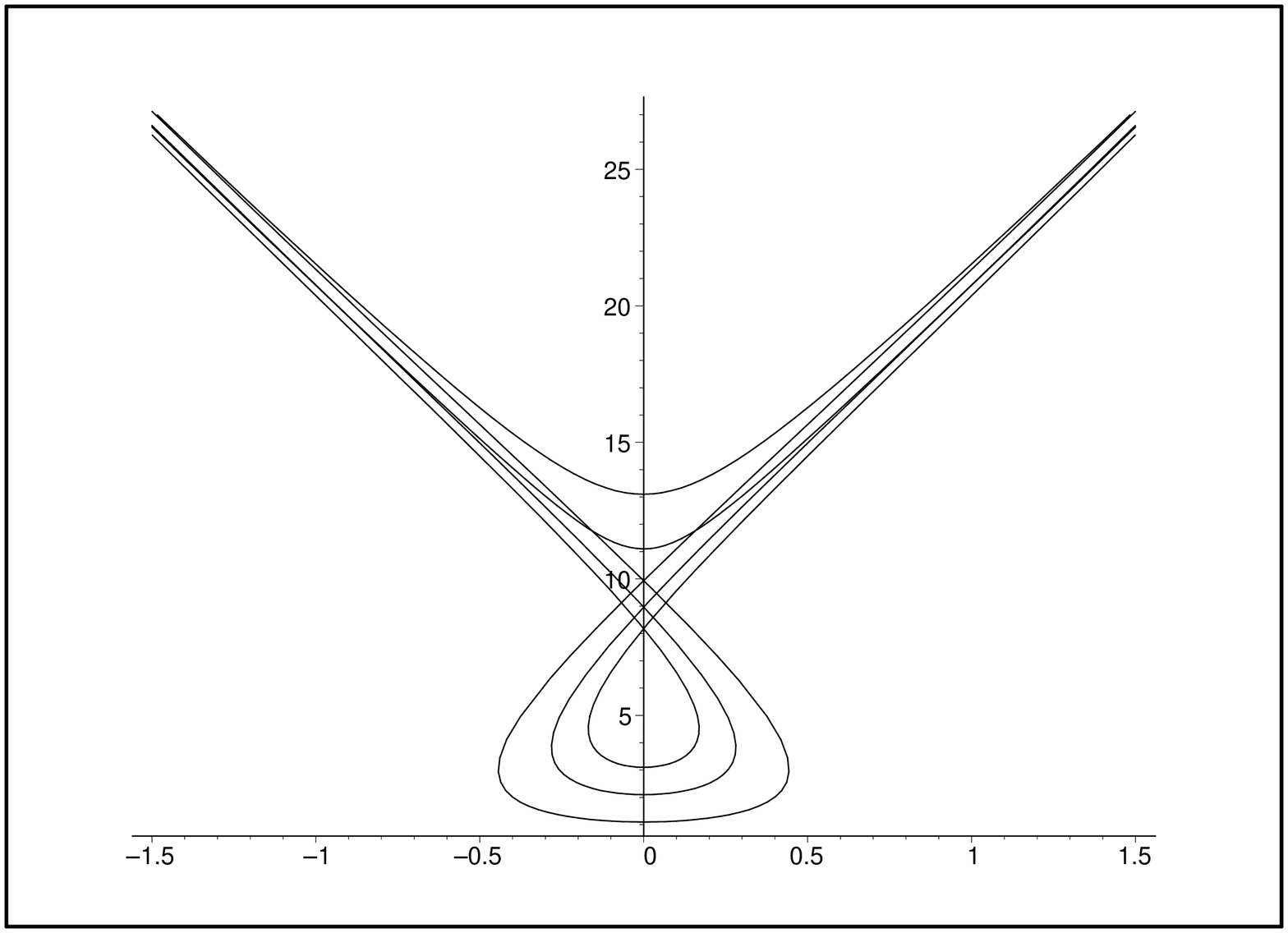}}}
\end{center}
\caption{A sketch of $V(x)$ for various values of $V_{0}$.}%
\label{V}%
\end{figure}

Equation (\ref{Eqtheta0}) admits a solution only when $0<V_{0}<2\pi$ (see
Figure \ref{Vo}). Thus, (\ref{EqV1}) generates a family of curves for $0<V_{0}<2\pi$
and it defines a single valued function $V(x)$ for $V_{0}%
\geq2\pi$ (see Figure \ref{V}). This was observed before in \cite{MR1376964}, but without any justification.

\section{Behavior for small $x$}

In this section we determine the expansion of $V(x)$ in a neighborhood of the
point $\left(  0,V_{0}\right)  .$ We shall distinguish two cases, depending on
$V_{0}\neq2\pi$ and $V_{0}=2\pi.$

\subsection{Case I: $V_{0}\neq2\pi$}

We consider the following expansion for $V(x)$%
\begin{equation}
V(x)=V_{0}+\sum\limits_{k=1}^{\infty}a_{k}\left(  \pi x\right)  ^{2k}%
\omega^{2k-1}\left(  -V_{0}\right)  ^{1-k}, \label{V1}%
\end{equation}
where%
\begin{equation}
\omega=\left[  \ln\left(  \frac{V_{0}}{2\pi}\right)  \right]  ^{-1}.
\label{omega}%
\end{equation}
Replacing (\ref{V1}) in (\ref{EqV1}) and expanding in powers of $x,$ we obtain%
\begin{align*}
a_{1}  &  =\omega\\
a_{2}  &  =\frac{4}{3}\omega^{2}%
\end{align*}%
\begin{align*}
a_{3}  &  =\frac{8}{15}\omega^{2}+\frac{28}{9}\omega^{3}\\
a_{4}  &  =\frac{32}{105}\omega^{2}+\frac{16}{5}\omega^{3}+\frac{80}{9}%
\omega^{4}%
\end{align*}%
\begin{align*}
a_{5}  &  =\frac{64}{315}\omega^{2}+\frac{528}{175}\omega^{3}+\frac{704}%
{45}\omega^{4}+\frac{2288}{81}\omega^{5}\\
a_{6}  &  =\frac{512}{3465}\omega^{2}+\frac{13312}{4725}\omega^{3}%
+\frac{14144}{675}\omega^{4}+\frac{5824}{81}\omega^{5}\\
&  +\frac{23296}{243}\omega^{6}%
\end{align*}%
\begin{align*}
a_{7}  &  =\frac{1024}{9009}\omega^{2}+\frac{12800}{4851}\omega^{3}%
+\frac{13312}{525}\omega^{4}+\frac{356864}{2835}\omega^{5}\\
&  +\frac{8704}{27}\omega^{6}+\frac{82688}{243}\omega^{7}\\
a_{8}  &  =\frac{4096}{45045}\omega^{2}+\frac{17408}{7007}\omega^{3}%
+\frac{5892608}{202125}\omega^{4}+\frac{1901824}{10125}\omega^{5}\\
&  +\frac{661504}{945}\omega^{6}+\frac{578816}{405}\omega^{7}+\frac
{909568}{729}\omega^{8}%
\end{align*}%
\begin{align*}
a_{9}  &  =\frac{8192}{109395}\omega^{2}+\frac{136192}{57915}\omega^{3}%
+\frac{153819136}{4729725}\omega^{4}+\frac{599911168}{2338875}\omega^{5}\\
&  +\frac{22685696}{18225}\omega^{6}+\frac{13536512}{3645}\omega^{7}%
+\frac{68913152}{10935}\omega^{8}+\frac{30764800}{6561}\omega^{9}\\
a_{10}  &  =\frac{131072}{2078505}\omega^{2}+\frac{3670016}{1640925}\omega
^{3}+\frac{229421056}{6449625}\omega^{4}+\frac{4575588352}{13820625}\omega
^{5}\\
&  +\frac{2095044608}{1063125}\omega^{6}+\frac{1036288}{135}\omega^{7}%
+\frac{208812032}{10935}\omega^{8}+\frac{6735872}{243}\omega^{9}\\
&  +\frac{117877760}{6561}\omega^{10}%
\end{align*}
and, in general, $a_{n}$ is a polynomial of degree $n$ in $\omega.$

\begin{figure}[ptb]
\begin{center}
\rotatebox{270} {\resizebox{12cm}{!}{\includegraphics{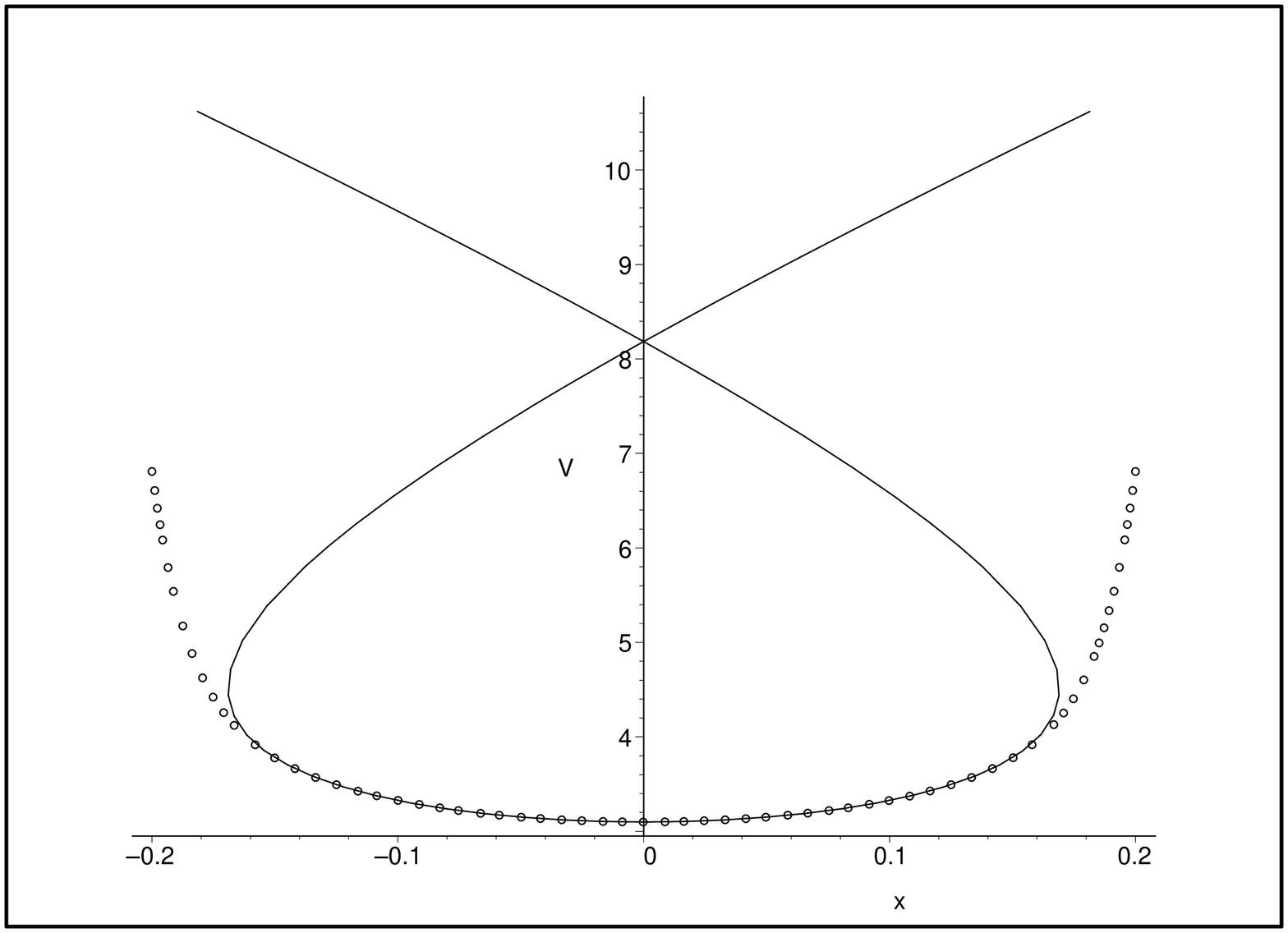}}}
\end{center}
\caption{A comparison of $V(x)$ (solid curve) and its approximation, with
$V_{0}=3.1$ (ooo).}%
\label{3}%
\end{figure}

\begin{figure}[ptb]
\begin{center}
\rotatebox{270} {\resizebox{12cm}{!}{\includegraphics{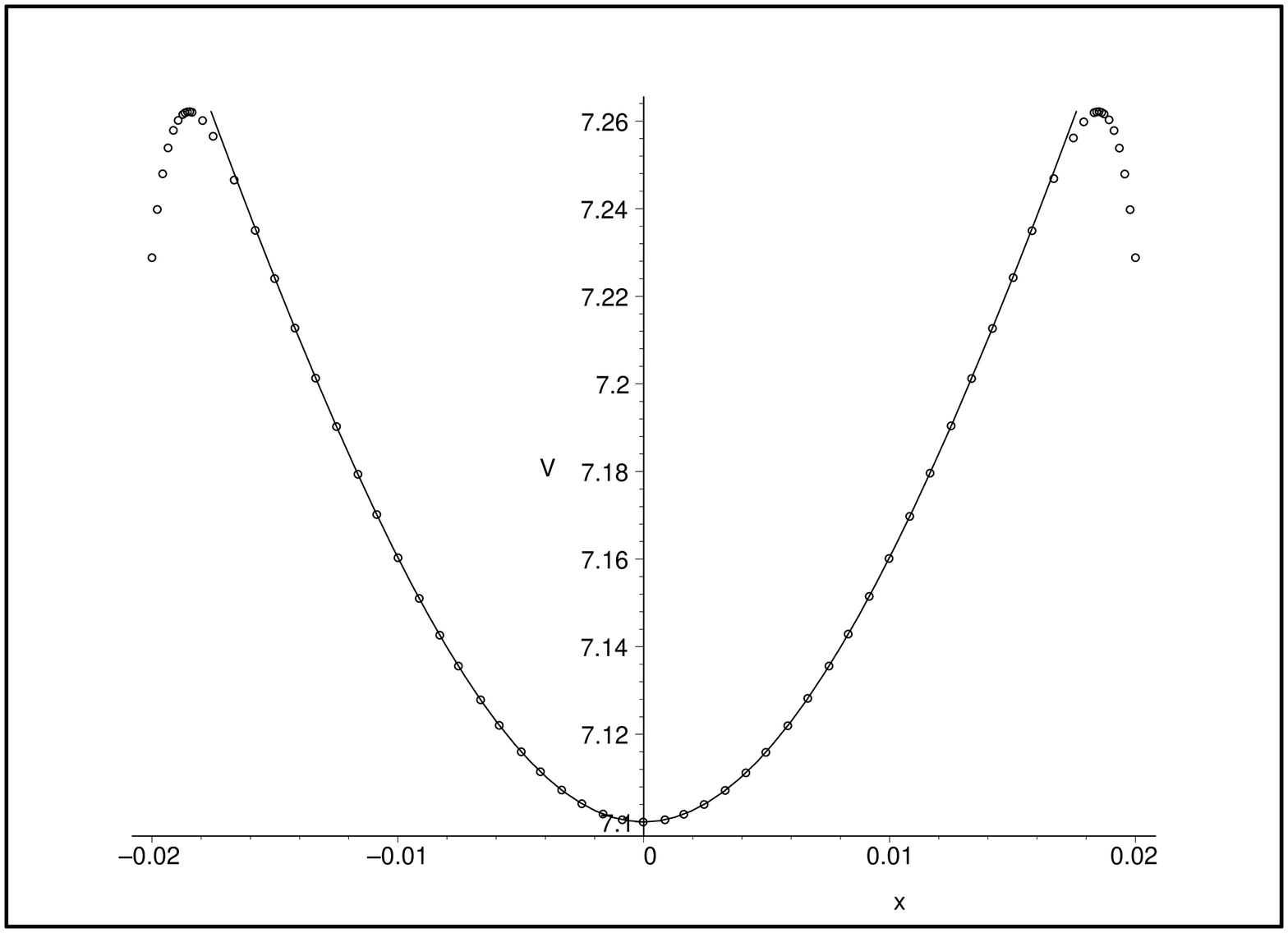}}}
\end{center}
\caption{A comparison of $V(x)$ (solid curve) and its approximation, with
$V_{0}=7.1$ (ooo).}%
\label{7}%
\end{figure}

In Figures \ref{3} and \ref{7}, we graph $V(x)$ and the first 10 terms of (\ref{V1}) with
$V_{0}=3.1$ and $V_{0}=7.1$, respectively.

\subsection{Case II: $V_{0}=2\pi$}

Since $\omega=\infty$ when $V_{0}=2\pi,$ the expansion (\ref{V1}) is no longer valid. Instead, we shall use%
\begin{equation}
V(x)=2\pi+\pi\sum\limits_{k=1}^{\infty}b_{k}\left(  3x\right)  ^{\frac{2k}{3}%
}5^{1-k}\pi^{\frac{k}{3}}. \label{V2}%
\end{equation}
Replacing (\ref{V2}) in (\ref{EqV1}) and expanding in powers of $x,$ we get%
\begin{align*}
b_{1}  &  =1,\quad b_{2}=\frac{2}{3},\quad b_{3}=\frac{1}{21},\quad
b_{4}=-\frac{2}{567},\quad b_{5}=\frac{92}{130977}\\
b_{6}  &  =-\frac{4}{21021},\quad b_{7}=\frac{19543}{321810489},\quad
b_{8}=-\frac{352610}{16412334939}\\
b_{9}  &  =\frac{12799}{1568439873},\quad b_{10}-\frac{6350075192}%
{1944910927276317},\ldots.
\end{align*}

\begin{figure}[ptb]
\begin{center}
\rotatebox{270} {\resizebox{12cm}{!}{\includegraphics{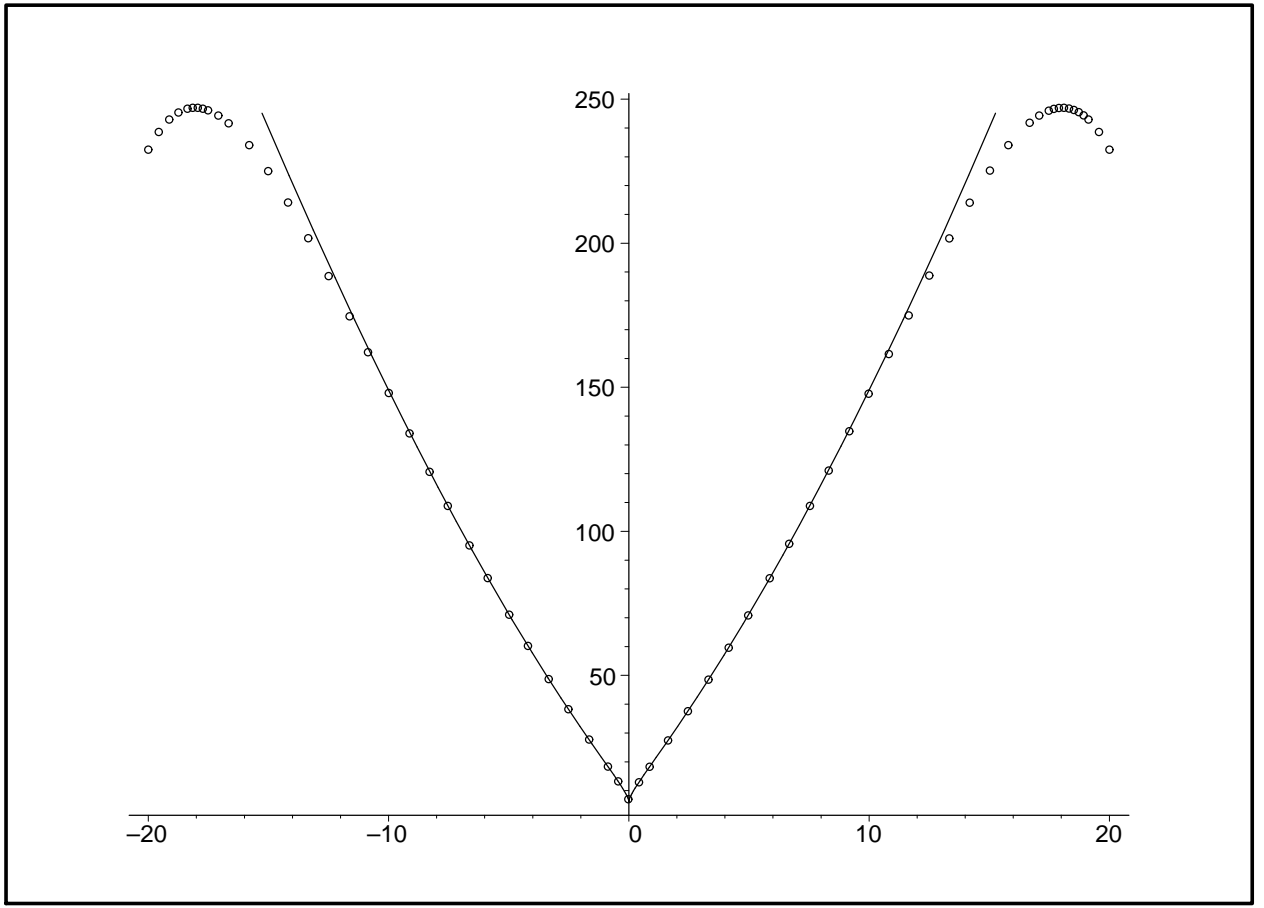}}}
\end{center}
\caption{A comparison of $V(x)$ (solid curve) and its approximation, with
$V_{0}=2\pi$ (ooo).}%
\label{2pi}%
\end{figure}

In Figure \ref{2pi} we sketch $V(x)$ and the the first 10 terms of (\ref{V2}).

\section{Asymptotic behavior}

We shall now analyze the behavior of $V(x)$ as $x\rightarrow\infty.$ From
(\ref{EqV1}) we have%
\[
\sqrt{V-V_{0}}\ln\left(  \frac{V_{0}}{2\pi e^{2}}\right)  +\sqrt{V}\ln\left[
\frac{\left(  \sqrt{V}+\sqrt{V-V_{0}}\right)  ^{2}}{V_{0}}\right]  \sim
\ln\left(  \frac{2V}{\pi e^{2}}\right)  \sqrt{V}%
\]
as $V\rightarrow\infty.$ Therefore,%
\[
\ln\left(  \frac{2V}{\pi e^{2}}\right)  \sqrt{V}\sim\pi x,\quad x\rightarrow
\infty,
\]
which implies%
\begin{equation}
V(x)\sim\frac{\pi^{2}x^{2}}{4}\left[  \mathrm{LW}\left(  \sqrt{\frac{\pi}{2}%
}\frac{\left\vert x\right\vert }{e}\right)  \right]  ^{-2},\quad\left\vert
x\right\vert \rightarrow\infty\label{Vasymp}%
\end{equation}
where $\mathrm{LW}\left(  \cdot\right)  $ represents the Lambert-W function
\cite{MR1414285}. Note that the asymptotic behavior is independent of $V_{0},$
as suggested by Figure \ref{V}.

\begin{figure}[ptb]
\begin{center}
\rotatebox{270} {\resizebox{12cm}{!}{\includegraphics{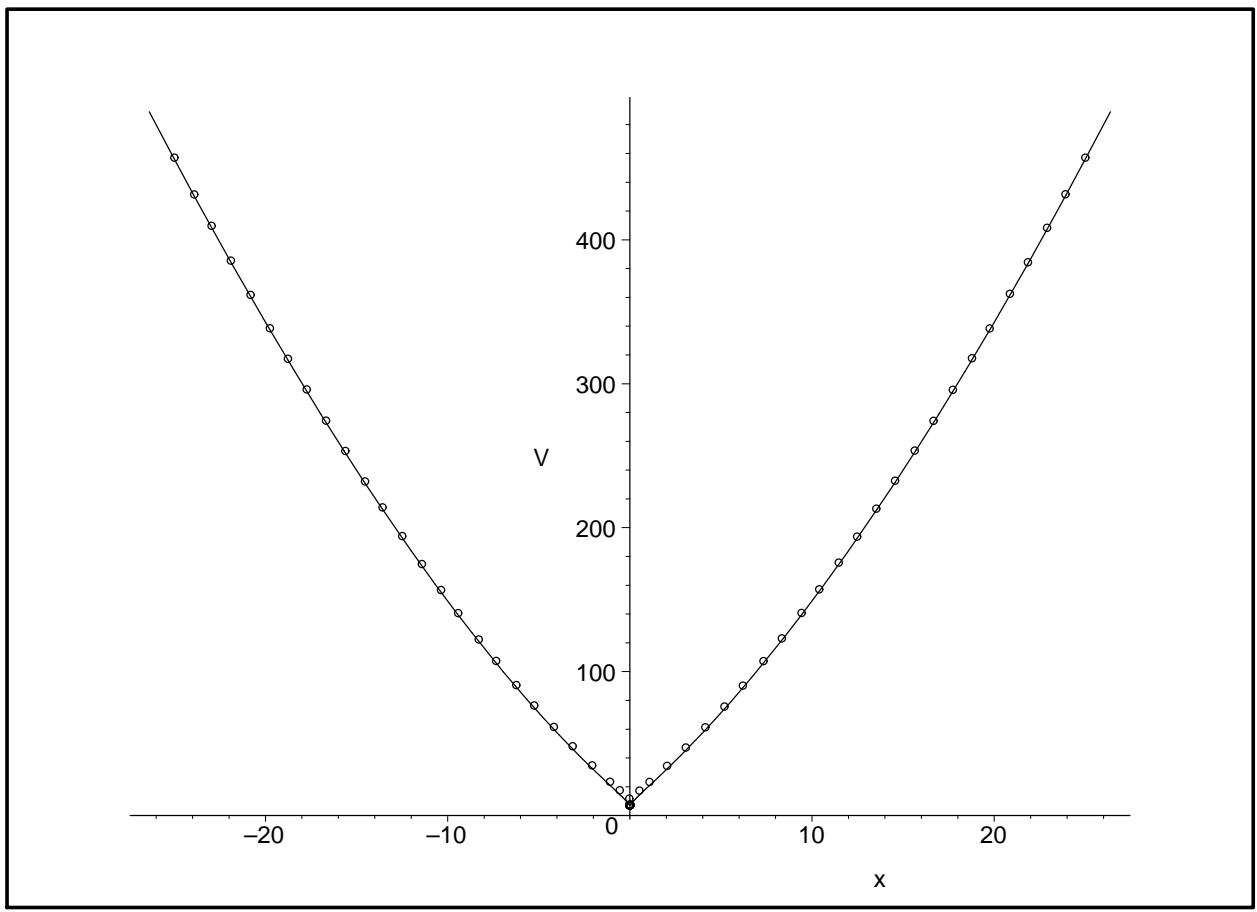}}}
\end{center}
\caption{A comparison of $V(x)$ (solid curve) and its asymptotic approximation (ooo).}%
\label{asym}%
\end{figure}

In Figure \ref{asym} we graph $V(x)$ and its asymptotic approximation (\ref{Vasymp}).

\end{document}